\documentclass[12pt]{report}
\usepackage{amssymb}
\usepackage{amsmath}

\begin{document}

\vspace{.2in}\parindent=0mm

\begin{center}

{\bf \large{VECTOR-VALUED NONUNIFORM MULTIRESOLUTION\\\parindent=0mm \vspace{.1in}  ASSOCIATED WITH LINEAR CANONICAL TRANSFORM }}

\vspace{.6in}\parindent=0mm

  {{M. YOUNUS BHAT}}

\parindent=0mm \vspace{.1in}
{{\it Department of  Mathematical Sciences,  Islamic University of Science and Technology Awantipora, Pulwama, Jammu and Kashmir 192122, India. E-mail: $\text{gyounusg@gmail.com}$}}

\parindent=0mm \vspace{.2in}
{{AAMIR H. DAR}}

\parindent=0mm \vspace{.1in}
{{\it Department of  Mathematical Sciences,  Islamic University of Science and Technology Awantipora, Pulwama, Jammu and Kashmir 192122, India. E-mail: $\text{gyounusg@gmail.com}$}}

\end {center}

\begin{quote}
\parindent=0mm \vspace{.2in}
 {A multiresolution analysis associated with linear canonical transform was defined by Shah and Waseem for which the translation set is a discrete set which is not a group. In this paper, we continue the study based on this nonstandard setting and introduce vector-valued nonuniform multiresolution analysis associated with linear canonical transform (LCT-VNUMRA) where the associated subspace $V_0^\mu$ of $L^2\big(\mathbb R, \mathbb C^{M}\big)$ has an orthonormal basis of the form $\left\{ {\bf \Phi} (x-\lambda)e^{-\frac{- \iota \pi A}{B}(t^2-\lambda^2)}\right\}_{ \lambda \in \Lambda}$  where $\Lambda = \left\{ 0,r/N\right\}+ 2 \mathbb Z, N \ge 1$ is an integer and $r$ is an odd integer such that $r$ and $N$ are relatively prime. We establish a necessary and sufficient condition for the existence of associated wavelets and derive an algorithm for the construction of vector-valued nonuniform multiresolution analysis on local fields starting from a vector refinement mask with appropriate conditions.}

\parindent=0mm\vspace{.1in}

\vspace{.1in}\parindent=0mm


\end{quote}

\parindent=0mm \vspace{.1in}
{\bf{1. Introduction}}

\parindent=0mm \vspace{.1in}
Multiresolution analysis (MRA) is an important mathematical tool since it provides a natural framework for understanding and constructing discrete wavelet systems. A multiresolution analysis is an increasing family of closed subspaces  $\left\{V_j: j\in\mathbb Z \right\}$ of $L^2(\mathbb R)$ such that $\bigcap_{j\in\mathbb Z}V_{j}=\left\{0\right\}, \, \bigcup_{j\in\mathbb Z}V_{j}$ is dense in $L^2(\mathbb R)$ and which satisfies $f\in V_j$ if and only if $f(2\cdot)\in V_{j+1}$. Furthermore, there exists an element $\varphi \in V_0$ such that the collection of integer translates of  function $\varphi,\,\left\{ \varphi(\cdot-k):k\in\mathbb Z\right\}$ represents a complete orthonormal system for $V_0$. The function $\varphi$ is called the {\it scaling function} or the {\it father wavelet}. The concept of multiresolution analysis has been extended in various ways in recent years. These concepts are generalized to  $L^2\big(\mathbb R^d\big)$, to lattices different from  $\mathbb Z^d$, allowing the subspaces of multiresolution analysis to be generated by Riesz basis instead of orthonormal basis, admitting a finite number of scaling functions, replacing the dilation factor 2 by an integer $M\geq 2$ or by an expansive matrix $A\in GL_{d}(\mathbb R)$ as long as $A\subset A\mathbb Z^d$. On the other hand, Xiang-Gen Xia and Suter$^{25}$ introduced the concept of vector-valued multiresolution analysis and orthogonal vector-valued wavelet basis and showed that vector-valued wavelets are a class of generalized multiwavelets. Chen and Cheng$^{4}$ presented the construction of a class of compactly supported orthogonal vector-valued wavelets and investigated the properties of vector-valued wavelet packets. Vector-valued wavelets are a class of generalized multiwavelets and multiwavelets can be generated from the component function in vector-valued wavelets. Vector-valued wavelets and multiwavelets are different in the following sense. Vector-valued wavelets can be used to decorrelate a vector-valued signal not only in the time domain but also between components for a fixed time where as multiwavelets focuses only on the decorrelation of signals in time domain. Moreover, prefiltering is usually required for discrete multiwavelet transform but not necessary for discrete vector-valued wavelet transforms. But all these concepts are developed on regular lattices, that is the translation set is always a group. Recently, Gabardo and Nashed$^{7,8}$ considered a generalization of Mallat's$^{13}$ celebrated theory of  multiresolution analysis based on spectral pairs, in which the translation set acting on the scaling function associated with the  multiresolution analysis to generate the subspace $V_0$ is no longer a group, but is the union of $\mathbb Z$ and a translate of $\mathbb Z$. More results in this direction can be found in Refs. 14, 22.

\pagestyle{myheadings}

\parindent=8mm \vspace{.2in}
The concept of novel multiresolution analysis in nonuniform settings was established by Shah and Waseem. They call it Nonuniform Multiresolution analysis associated with linear canonical transform (LCT-NUMRA). They also constructed associated wavelet packets and presented orthogonal decomposition. In this paper, we continue the study based on this nonstandard setting and introduce vector-valued nonuniform multiresolution analysis associated with linear canonical transform (LCT-VNUMRA) where the associated subspace $V_0^\mu$ of $L^2\big(\mathbb R, \mathbb C^{M}\big)$ has an orthonormal basis of the form $\left\{ {\bf \Phi} (x-\lambda)e^{-\frac{- \iota \pi A}{B}(t^2-\lambda^2)}\right\}_{ \lambda \in \Lambda}$  where $\Lambda = \left\{ 0,r/N\right\}+ 2 \mathbb Z, N \ge 1$ is an integer and $r$ is an odd integer such that $r$ and $N$ are relatively prime. We establish a necessary and sufficient condition for the existence of associated wavelets and derive an algorithm for the construction of vector-valued nonuniform multiresolution analysis on local fields starting from a vector refinement mask with appropriate conditions.

\parindent=8mm \vspace{.2in}
This paper is organized as follows.  In Sec. 3, we review the uniform and non-uniform multiresolution analysis associated with LCT and  certain properties related to the construction of associated wavelets. In Sec. 4, we introduce the notion of vector-valued nonuniform multiresolution analysis associated with linear canonical transform (LCT-VNUMRA) and establish a necessary and sufficient condition  for the existence of associated wavelet. In Sec. 5, we construct a  LCT-VNUMRA starting from a vector refinement mask  satisfying appropriate conditions.

\parindent=0mm \vspace{.2in}
{\bf{3. Nonuniform Multiresolution Analysis Associated with Linear Canonical Transform}}

\parindent=8mm \vspace{.2in}
 For the sake of simplicity, we consider the second order matrix $\mu_{2\times 2}=(A, B, C, D)$ with its transpose defined by  $\mu_{2\times 2}^T=(A, B, C, D)^T$. Let us first introduce the definition of Linear Canonical Transform.
 
\parindent=0mm \vspace{.2in}
{\bf{Definition 3.1.}} The linear canonical transform of any $ f \in L^2(\mathbb R)$ with respect to the unimodular matrix $\mu_{2\times 2}=(A, B, C, D)$is defined by
 $${\mathcal L}[f](\xi)=\left\{\begin{array}{lcr}
 \int_\mathbb R f(t){\mathcal K}_\mu(t, \xi) dt&&B \ne 0\\\\
 \sqrt D \exp {\frac{CD \xi^2}{2}}f(D\xi)&& B=0.
 \end{array}\right.$$
 where ${\mathcal K}_\mu(t, \xi)$ is the kernel of linear canonical transform and is given by
 $${\mathcal K}_\mu(t, \xi)=\dfrac{1}{\sqrt{2 \pi \iota B}}\exp \left\{\dfrac{\iota(At^2-2 t \xi+D \xi^2)}{2B}\right\}, \quad B \ne 0$$
\parindent=8mm \vspace{.2in}
Recently, Shah and Waseem$^{22}$ considered a generalization of the notion of multiresolution analysis associated with linear canonical transform, which is called {\it nonuniform multiresolution analysis associated with linear canonical transform } (LCT-NUMRA) and is based on the theory of spectral pairs. In this set up, the associated subspace $V_0^\mu$ of $L^2(\mathbb R)$ has an orthonormal basis, a collection of translates of the scaling function $\varphi$ of the form $\left\{ \varphi (t-\lambda)e^{-\frac{- \iota \pi A}{B}(t^2-\lambda^2)}\right\}_{ \lambda \in \Lambda}$  where $\Lambda = \left\{ 0,r/N\right\}+ 2 \mathbb Z, N \ge 1$ is an integer and $r$ is an odd integer such that $r$ and $N$ are relatively prime.

\parindent=8mm \vspace{.1in}
We first recall the definition of a nonuniform multiresolution analysis associated with linear canonical transform (as defined in Ref. 22) and
some of its properties.

\parindent=0mm \vspace{.2in}
{\bf{Definition 3.2.}} For an integer $N \ge 1$ and an odd integer $r$ with $1\leq r \leq 2N-1$ such that $r$ and $N$ are relatively prime, a  nonuniform multiresolution analysis associated with linear canonical transform is a sequence of closed subspaces $\left\{V_j^\mu: j\in\mathbb Z\right\}$ of $L^2(\mathbb R)$ such that the following properties hold:

\parindent=0mm \vspace{.2in}
(a)\quad $V_j^\mu \subset V_{j+1}^\mu\; \text{for all}\; j \in \mathbb Z;$

\parindent=0mm \vspace{.1in}
(b)\quad $\bigcup_{j\in \mathbb Z}V_j^\mu\;\text{is dense in}\;L^2(K);$

\parindent=0mm \vspace{.1in}
(c)\quad $\bigcap_{j\in \mathbb Z}V_j^\mu=\{0\};$

\parindent=0mm \vspace{.1in}
(d)\quad $f(t) \in V_j^\mu\; \text{if and only if}\;f\left(2N\cdot\right)e^{-\iota \pi A\left(1-(2N)^2\right)t^2/B} \in V_{j+1}^\mu\; \text{for all}\; j \in \mathbb Z;$

\parindent=0mm\vspace{.1in}

(e)~ There exists a function $\varphi$ in $V_0^\mu$ such that $\left\{ \varphi (t- \lambda )e^{-\frac{- \iota \pi A}{B}(t^2-\lambda^2)}: \lambda \in \Lambda\right\}$, is a complete orthonormal basis for $V_0^\mu$.

\parindent=8mm \vspace{.2in}
Given a LCT- NUMRA $\left\{V_{j}^\mu: j\in\mathbb Z\right\}$, we define another sequence $\left\{W_{j}^\mu: j\in\mathbb Z\right\}$ of closed subspaces of $L^2(\mathbb R)$ by $W_{j}^\mu:=V_{j+1}^\mu\ominus V_{j}^\mu, j\in\mathbb Z.$ These subspaces inherit the scaling property of $V_{j}^\mu$, namely,

$$f(\cdot)\in W_j^\mu\quad \text{if and only if } \quad f\left(2N\cdot\right)e^{2\iota \pi \lambda \xi/B}\in W_{j+1}^\mu.\eqno (3.1)$$

\parindent=0mm \vspace{.1in}
Moreover, the subspaces $\left\{W_{j}^\mu: j\in\mathbb Z\right\}$ are mutually orthogonal, and we have the following orthogonal decomposition:

$$L^2(R)=\bigoplus_{j\in\mathbb Z}W_{j}^\mu=V_{0}^\mu \oplus\left(\bigoplus_{j\ge 0}W_{j}^\mu \right).\eqno (3.2)$$

\parindent=8mm \vspace{.1in}
 A set of functions $\left\{\psi_{1}^\mu, \psi_{1}^\mu,\dots,\psi_{2N-1}^\mu\right\}$ in $L^2(K\mathbb R$ is said to be a {\it set of  basic wavelets} associated with the LCT-NUMRA $\left\{V_{j}^\mu: j\in\mathbb Z\right\}$ if the family of functions $\left\{ \psi_{\ell}(t-\lambda)e^{-\frac{- \iota \pi A}{B}(t^2-\lambda^2)}:1\le \ell\le 2N-1, \lambda\in \Lambda\right\}$ forms an orthonormal basis for $W_{0}^\mu$.

\parindent=8mm \vspace{.1in}
In view of (3.1) and (3.2), it is clear that if $\left\{\psi_{1}, \psi_{1},\dots,\psi_{2N-1}\right\}$ is a set of basic wavelets, then $\left\{ (2N)^{j/2}\psi_{\ell}\big((2N)^{j}t-\lambda\big)e^{-\frac{- \iota \pi A}{B}(t^2-\lambda^2)}:1\le \ell\le 2N-1, \lambda\in \Lambda\right\}$ constitutes an orthonormal basis for $L^2(K)$.

\parindent=0mm \vspace{.2in}
{\bf{4. Vector-valued Nonuniform Multiresolution Associated with Linear Canonical Transform }}

\parindent=0mm \vspace{.2in}
In this section, we introduce the notion of vector-valued nonuniform multiresolution analysis associated with linear canonical transform and establish a necessary and sufficient condition for the existence of associated wavelets.

\parindent=8mm \vspace{.1in}
Let $M$ be a constant and $2\le M\in\mathbb Z$. By $L^2\big(\mathbb R,\mathbb C^M\big)$, we denote the set of all vector-valued functions ${\bf f}(x)$ i.e.,
$$L^2\big(\mathbb R,\mathbb C^M\big)= \left\{ {\bf f}(x)= \big( f_{1}(x),f_{2}(x),\dots,f_{M}(x)\big)^{T} :x\in \mathbb R, f_{t}(x)\in L^2(\mathbb R),t=1,2,\dots, M \right\}$$

\parindent=0mm \vspace{.1in}
where $T$ means the transpose of a vector. The space $L^2\big(\mathbb R,\mathbb C^M\big)$ is called {\it vector-valued function space}. For ${\bf f}(x)\in L^2\big(\mathbb R,\mathbb C^M\big),\, \big\|{\bf f}\big\|$ denotes the norm of vector-valued function ${\bf f}$ and is defined as:

$$\big\|{\bf f}\big\|_{2}= \left(\sum_{t=1}^{M} \int_{\mathbb R} \big|f_{t}(x)\big|^2 dx\right)^{1/2}.\eqno(4.1)$$

\parindent=0mm \vspace{.1in}
For a vector-valued function ${\bf f}(x)\in L^2\big(\mathbb R,\mathbb C^M\big)$, the integration of ${\bf f}(x)$ is defined as:

$$\int_{\mathbb R}{\bf f}(x) dx =\left( \int_{\mathbb R}f_{1}(x) dx, \int_{\mathbb R}f_{2}(x) dx, \dots, \int_{\mathbb R}f_{M}(x) dx \right)^{T}.$$

\parindent=0mm \vspace{.1in}
For any two vector-valued functions ${\bf f, g}\in L^2\big(\mathbb R,\mathbb C^M\big)$, their vector-valued inner product $\langle {\bf f}, {\bf g}\rangle$ is defined as:

$$\langle {\bf f}, {\bf g}\rangle=\int_{\mathbb R} {\bf f}(x)\overline{{\bf g}(x)}\, dx.\eqno(4.2)$$

\parindent=0mm \vspace{.1in}
With $\Lambda = \left\{ 0, r/N\right\}+2 \mathbb Z$ as defined above, we  define the {\it vector-valued nonuniform multiresolution analysis associated with linear canonical transform } (LCT-VNUMRA) as follows:

\parindent=0mm \vspace{.2in}
{\bf{Definition 4.1.}}Given a real uni-modular matrix $\mu=(A, B, C, D)$ and an integer $N \ge 1$ and an odd integer $r$ with $1\leq r \leq 2N-1$ such that $r$ and $N$ are relatively prime, an associated linear canonical  vector-valued non-uniform multiresolution analysis (LCT-LCT-VNUMRA) is a sequence of closed subspaces $\left\{V_j^\mu: j\in\mathbb Z\right\}$ of $L^2\big(\mathbb R,\mathbb C^M\big)$ such that the following properties hold:

\parindent=0mm \vspace{.2in}
(a)\quad $V_j^\mu \subset V_{j+1}^\mu\; \text{for all}\; j \in \mathbb Z;$

\parindent=0mm \vspace{.1in}
(b)\quad $\bigcup_{j\in \mathbb Z}V_j^\mu\;\text{is dense in}\;L^2\big(\mathbb R,\mathbb C^M\big);$

\parindent=0mm \vspace{.1in}
(c)\quad $\bigcap_{j\in \mathbb Z}V_j^\mu=\{\bf 0\}$,\; where ${\bf 0}$ is the zero vector of $L^2\big(\mathbb R,\mathbb C^M\big);$

\parindent=0mm \vspace{.1in}
(d)\quad ${\bf \Phi}(t) \in V_j^\mu\; \text{if and only if}\;{\bf \Phi}(2N t)e^{-\iota \pi A\left(1-(2N)^2\right)t^2/B} \in V_{j+1}^\mu\; \text{for all}\; j \in \mathbb Z;$

\parindent=0mm\vspace{.1in}

(e)~ There exists a function $\Phi$ in $V_0^\mu$ such that $\left\{\Phi_{0, \lambda}^\mu(t)= {\bf \Phi} (t- \lambda )e^{-\frac{- \iota \pi A}{B}(t^2-\lambda^2)}: \lambda \in \Lambda\right\}$, is a complete orthonormal basis for $V_0^\mu$. The vector valued function ${\bf \Phi}(x)$ is called a {\it vector-valued scaling function} of the LCT-VNUMRA.

\parindent=8mm \vspace{.2in}
For every $j\in\mathbb Z$, define $W_{j}^\mu$ to be the orthogonal complement of $V_{j}^\mu$ in $V_{j+1}^\mu$. Then we have
$$V_{j+1}^\mu=V_{j}^\mu \oplus W_{j}^\mu \quad \text{and}\quad W_{\ell}^\mu\perp W_{\ell}^{\prime \mu}\quad \text{if}~\ell\ne \ell^\prime.\eqno(4.3) $$

\parindent=0mm \vspace{.1in}
It follows that for $j>J$,
$$V_{j}^\mu=V_{J}^\mu\oplus \bigoplus_{\ell=0}^{j-J-1}W_{j-\ell}^\mu\,,\eqno(4.3) $$

\parindent=0mm \vspace{.1in}
where all these subspaces are orthogonal. By virtue of condition (b) in the Definition 4.1, this implies
$$L^2\big(\mathbb R,\mathbb C^M\big)=\bigoplus_{j\in\mathbb Z}W_{j}^\mu,\eqno(4.5) $$

\parindent=0mm \vspace{.1in}
a decomposition of $L^2\big(\mathbb R,\mathbb C^M\big)$ into mutually orthogonal subspaces.

\parindent=8mm \vspace{.2in}
As in the standard case, one expects the existence of $2N -1$ number of functions so that their translation by elements of $\Lambda$ and dilations by the integral powers of ${\mathfrak p^{-1}}N$ form an orthonormal basis for $L^2\big(\mathbb R,\mathbb C^M\big)$.

\parindent=0mm \vspace{.2in}
{\bf{Definition 4.2.}} A set of functions $\left\{{\bf \Psi}_{1}^\mu, {\bf \Psi}_{2}^\mu,\dots,{\bf \Psi}_{2N-1}^\mu\right\}$ in $L^2\big(\mathbb R,\mathbb C^M\big)$ will be called a {\it set of basic wavelets} associated with a given LCT-VNUMRA if the family of functions $\left\{ {\bf \Psi}_{\ell}(t-\lambda)e^{-\frac{- \iota \pi A}{B}(t^2-\lambda^2)}:1\le \ell\le 2N-1, \lambda\in \Lambda\right\}$ forms an orthonormal basis for $W_{0}^\mu$.

\parindent=5mm \vspace{.2in}
In the following, we want to seek a set of wavelet functions $\left\{{\bf \Psi}_{1}^\mu, {\bf \Psi}_{2}^\mu,\dots,{\bf \Psi}_{2N-1}^\mu\right\}$ in $ W_{0}^\mu$ such that $\left\{ (2N)^{j/2}{\bf \Psi}_{\ell}\big( (2N)^{j}t-\lambda\big)e^{-\frac{- \iota \pi A}{B}(t^2-\lambda^2)}:1\le \ell\le 2N-1, \lambda\in \Lambda\right\}$ form an orthonormal basis of $W_{j}^\mu$. By the nested structure of LCT-LCT-VNUMRA, this task can be reduce to find ${\bf \Psi}_{\ell}\mu\in W_{0}^\mu$ such that $\left\{{\bf \Psi}_{\ell}\big(t-\lambda\big)e^{-\frac{- \iota \pi A}{B}(t^2-\lambda^2)}:1\le \ell\le 2N-1, \lambda\in \Lambda\right\}$ constitutes an orthonormal basis of $W_{0}^\mu$.

\parindent=8mm \vspace{.2in}
Let ${\bf \Phi}=\big(\varphi_1^\mu, \varphi_2^\mu,\dots,\varphi_M^\mu\big)^T$ be a scaling vector of the given LCT-VNUMRA. Since $
{\bf \Phi}\in V_{0}^\mu\subset V_{1}^\mu$, there exist  $M\times M$ constant matrix sequence $\{G_\lambda\}_{\lambda \in \Lambda}$ such that

$${\bf \Phi}(t)=\sqrt {2N}\sum_{\lambda \in \Lambda}G_\lambda{\bf \Phi}\big(2Nt-\lambda \big)e^{-\frac{- \iota \pi A}{B}(t^2-\lambda^2)}. \eqno(4.6)$$
where $G_\lambda=\displaystyle\int_{\mathbb R}{\bf \Phi}(t)e^{-\iota \pi A\left(1-(2N)^2\right)t^2/B}\overline{\phi_{1,\lambda}^\mu(t)}dt $.

\parindent=0mm \vspace{.2in}
Taking linear canonical transform on both sides of equation (4.6), we obtain

$$\mathcal L[{\bf \Phi}(t)](\xi)={\hat {\bf \Phi}}\left(\frac{\xi}{B}\right)=G^\mu\left( \frac{\xi}{2NB}\right){\hat {\bf \Phi}}\left( \frac{\xi}{2NB}\right),\eqno(4.7)$$

\parindent=0mm \vspace{.1in}
where $G^\mu\left(\frac{\xi}{B}\right)=\frac{1}{\sqrt {2N}}\sum_{\lambda \in \Lambda}G_\lambda^\mu\,\overline{e^{-2\pi \iota \lambda \xi /B}},$  is called {\it symbol} or {\it vector refinement mask} of the scaling function $\Phi$. By replacing $\xi$ by ${\xi}/{2NB}$ in relation (4.7), we obtain

$$\hat {\bf \Phi}\left( \frac{\xi}{2NB}\right)=G^\mu\left( \left(\frac{1}{2NB}\right)^2\xi\right)\hat {\bf \Phi}\left( \left(\frac{1}{2NB}\right)^2\xi\right),$$

\parindent=0mm \vspace{.0in}
and then

$$\hat {\bf \Phi}(\xi)=G^\mu\left( \frac{\xi}{2NB}\right)G^\mu\left( \left(\frac{1}{2NB}\right)^2\xi\right)\hat {\bf \Phi}\left( \left(\frac{1}{2NB}\right)^2\xi\right).$$

\parindent=0mm \vspace{.1in}
We can continue this and obtain, for any $n\in\mathbb N$,

$$\begin{array}{rcl}
\hat {\bf \Phi}(\xi)&=&\displaystyle G^\mu\left( \frac{\xi}{2NB}\right)G^\mu\left( \left(\frac{1}{2NB}\right)^2\xi\right)\cdots G^\mu\left( \left(\frac{1}{2NB}\right)^n\xi\right)\hat {\bf \Phi}\left( \left(\frac{1}{2NB}\right)^n\xi\right)\\\\
&=&\displaystyle \hat {\bf \Phi}\left( \left(\frac{1}{2NB}\right)^n\xi\right) \prod_{m=1}^{n}G^\mu\left( \left(\frac{1}{2NB}\right)^m\xi\right).
\end{array}$$

\parindent=0mm \vspace{.1in}
By taking  $n\to \infty$  and noting that $\left|\left(\dfrac{1}{2NB}\right)^n\right|= \dfrac{1}{(2NB)^n}\to 0$ as $n\to \infty$, the above relation reduces to
$$\hat {\bf \Phi}(\xi)= \hat {\bf \Phi}(0) \prod_{m=1}^{\infty} G^\mu \left( \left(\frac{1}{2NB}\right)^m\xi\right).\eqno(4.8)$$

\parindent=0mm \vspace{.1in}
As usual, we assume $\hat {\bf \Phi}(\xi)$ is continuous at zero, and $\hat{\bf \Phi}(0)=I_{M}$, where $I_{M}$ denotes the identity matrix of order $M\times M$. Therefore, equation (4.8) becomes

$$\hat{\bf \Phi}(\xi)= \prod_{m=1}^{\infty} G^\mu\left( \left(\frac{1}{2NB}\right)^m\xi\right)\eqno(4.9)$$

\parindent=0mm \vspace{.1in}
Moreover, it is immediate from (4.7) that $G(0)=I_{M}$, which is essential for convergence of the infinite product $\prod_{m=1}^{\infty} G^\mu \left( \left(\frac{1}{2NB}\right)^m\xi\right)$.

\parindent=8mm \vspace{.2in}
We now investigate the orthogonal property of the scaling function ${\bf \Phi}$ by means of the vector refinement mask $G(\xi)$.

\parindent=0mm \vspace{.2in}
{\bf{Lemma 4.3.}} {\it If ${\bf \Phi} \in L^2\big(\mathbb R,\mathbb C^M\big)$ defined by Equation (4.6) is an orthogonal vector-valued scaling function, then we have

$$\sum_{m \in 2\mathbb Z}G^\mu_{m}\overline{G^\mu_{2NB(\lambda-\lambda^\prime)+m}}=2NB\delta_{\lambda, \lambda^\prime}I_M,\quad \forall~\lambda, \lambda^\prime\in \Lambda, \eqno(4.10)$$

\parindent=0mm \vspace{.0in}
where $\delta_{\lambda, \lambda^\prime}$ denotes the Kronecker's delta.}

\parindent=0mm \vspace{.2in}
{\bf{Proof.}} Since the scaling function is orthogonal vector-valued, we have
$$\begin{array}{rcl}
\delta_{\lambda, \lambda^\prime}I_M &=&\displaystyle\int_\mathbb R{\bf \Phi}(t-\lambda)e^{-\frac{- \iota \pi A}{B}(t^2-\lambda^2)}\overline{{\bf \Phi}(t-\lambda^\prime)e^{-\frac{- \iota \pi A}{B}(t^2-\lambda{^ \prime 2})}}dt\\\\
&=&\displaystyle\sum_{\sigma \in \Lambda}\displaystyle\int_\mathbb R G^\mu_\sigma{\bf \Phi}\big(2NBt-2NB\lambda-\sigma\big)\displaystyle\sum_{\sigma \in \Lambda}\overline{G^\mu_{\sigma }}\,\overline{{\bf \Phi}\big(2NBt-2NB\lambda^\prime-\sigma \big)}dt\\\\
&=&\displaystyle\sum_{\sigma \in \Lambda}\displaystyle\sum_{\sigma  \in \Lambda}G^\mu_\sigma\left\{\displaystyle\int_\mathbb R{\bf \Phi}\big(2NBt-2NB\lambda-\sigma\big)\overline{{\bf \Phi}\big(2NBt-2NB\lambda^\prime-\sigma \big)}\right\}dt\overline{G^\mu_{\sigma }}\\\\
&=&\dfrac{1}{2NB}\displaystyle\sum_{\sigma \in \Lambda}\displaystyle\sum_{\sigma  \in \Lambda}G^\mu_\sigma\left\{\displaystyle\int_\mathbb R{\bf \Phi}\big(t-2NB\lambda-\sigma\big)\overline{{\bf \Phi}\big(t-2NB\lambda^\prime-\sigma \big)}\right\}dt \overline{G^\mu_{\sigma }}.
\end{array}$$

\parindent=0mm \vspace{.1in}
Taking $\sigma=2m$ and $\sigma=2n$, where $m,n \in \mathbb Z$, we have
$$\begin{array}{rcl}
\delta_{\lambda, \lambda^\prime}I_M&=&\dfrac{1}{2NB}\displaystyle\sum_{\sigma \in \Lambda}\displaystyle\sum_{\sigma  \in \Lambda}G^\mu_\sigma\left\langle{\bf \Phi}\big(t-2NB\lambda-\sigma\big),\overline{{\bf \Phi}\big(t-2NB\lambda^\prime-\sigma \big)}\right\rangle \overline{G^\mu_{\sigma }}\\\\
&=&\dfrac{1}{2NB}\displaystyle\sum_{m \in \mathbb N_0}\displaystyle\sum_{n \in \mathbb N_0}G^\mu_{2m}\left\langle{\bf \Phi}\big(t-2NB\lambda-2m\big),\overline{{\bf \Phi}\big(t-2NB\lambda^\prime-2n\big)}\right\rangle \overline{G^\mu_{2n }}\\\\
&=&\dfrac{1}{2NB}\displaystyle\sum_{m \in \mathbb N_0}G^\mu_{2m}\overline{G^\mu_{2NB(\lambda-\lambda^\prime)+2m}}.\end{array}$$

\parindent=0mm \vspace{.1in}
Therefore, identity (4.10) follows.

\parindent=0mm \vspace{.2in}
Taking $\sigma=\dfrac{r}{NB}+2m$ and $\sigma=2n$, where $m,n \in \mathbb Z$, we have
$$\begin{array}{rcl}
\delta_{\lambda, \lambda^\prime}I_M&=&\dfrac{1}{2NB}\displaystyle\sum_{\sigma \in \Lambda}\displaystyle\sum_{\sigma  \in \Lambda}G^\mu_\sigma\left\langle{\bf \Phi}\big(t-2NB\lambda-\sigma\big),\overline{{\bf \Phi}\big(t-2NB\lambda^\prime-\sigma \big)}\right\rangle \overline{G^\mu_{\sigma }}\\\\
&=&\dfrac{1}{2NB}\displaystyle\sum_{m \in 2\mathbb Z}\displaystyle\sum_{n \in 2\mathbb Z}G^\mu_{\frac{r}{NB}+2m}\Big\langle{\bf \Phi}\big(t-2NB\lambda-\frac{r}{NB}-2m\big),\\\
&&\qquad\qquad\qquad\qquad\qquad\qquad\qquad\overline{{\bf \Phi}\big(t-2NB\lambda^\prime-2n\big)}\Big\rangle \overline{G^\mu_{2n }}\\\
&=&\dfrac{1}{2NB}\displaystyle\sum_{m \in 2\mathbb Z}G^\mu_{2m}\overline{G^\mu_{2NB(\lambda-\lambda^\prime)+2m }}.\end{array}$$

\parindent=0mm \vspace{.1in}
Thus, in both the cases, we get the desired result. \qquad\qquad \fbox

\parindent=8mm \vspace{.2in}

We denote ${\bf \Psi}_{0}={\bf \Phi}$, the scaling function, and consider $2N-1$ functions ${\bf \Psi}^\mu_{\ell}, 1\le \ell\le 2N-1$, in $W^\mu_{0}$ as possible candidates for wavelets. Since $(1/2NB)\,{\bf \Psi}^\mu_{\ell}(1/2NB t)\in V^\mu_{-1}\subset V^\mu_{0}$, it follows from property (d) of Definition 4.1 that for each $\ell, ~0\le \ell\le 2N-1$, there exists a uniquely supported sequence $\left\{H_\lambda^\ell\right\}_{\lambda \in \Lambda,\,1\le \ell \le 2N-1}$ of $M\times M$ constant matrices such that

$${\bf \Psi}^\mu_\ell(t)=\sqrt {2N} \sum_{\lambda \in \Lambda}H_{\lambda,\ell}\,{\bf \Phi}\left(2Nt-\lambda\right)e^{-\frac{- \iota \pi A}{B}(t^2-\lambda^2)}. \eqno(4.11)$$

\parindent=0mm \vspace{.1in}
On taking the linear canonical  transform on both sides of Equation (4.11), we have

$$\hat{\bf \Psi}^\mu_\ell\left(2N\frac{\xi}{B}\right)=H^\mu_\ell\left(\frac{\xi}{B}\right)\hat{\bf \Phi}\left(\frac{\xi}{B}\right),\eqno(4.12)$$
where
$$H^\mu_\ell\left(\frac{\xi}{B}\right)=\dfrac{1}{\sqrt {2N}} \sum_{\lambda \in \Lambda}H_{\lambda, \ell}^\mu\,e^{-2 \pi \iota \lambda \xi /B}.\eqno(4.13)$$

\parindent=0mm \vspace{.1in}
In view of the specific form of $\Lambda = \left\{ 0, \dfrac{r}{N}\right\}+{2 \mathbb Z}$,, we observe that

$$H^\mu_\ell\left(\frac{\xi}{B}\right)=H^{\mu, 1}_\ell\left(\frac{\xi}{B}\right)+e^{-2 \pi \iota  \xi /NB}H^{\mu, 2}_\ell\left(\frac{\xi}{B}\right),\quad  0\le \ell\le 2N-1,\eqno(4.14) $$

\parindent=0mm \vspace{.1in}
where $H^{\mu, 1}_\ell$ and $H^{\mu, 2}_\ell$ are $M \times M$ constant symmetric matrix sequences.

\parindent=0mm \vspace{.2in}
{\bf{Lemma 4.4.}} {\it Consider a LCT-VNUMRA on  $\mathbb R$ as in Definition 4.1. Suppose that there exist $2 N-1$ functions ${\bf \Psi}_{k},\, k=1,2,\dots, 2N-1$ in $V_1$. Then the family of functions

$$\big\{ {\bf \Psi}_k(t- \lambda )e^{-\frac{- \iota \pi A}{B}(t^2-\lambda^2)}:\lambda \in \Lambda,\,k=0, 1, \dots, 2 N-1  \big\}\eqno(4.15)$$

\parindent=0mm \vspace{.1in}
 forms an orthonormal system in $V_1$ if and only if}

$$\sum_{r=0}^{2N-1}H^\mu_k\left( \frac{\xi}{2NB} +\frac{r}{4N}\right)\overline{H^\mu_\ell\left( \frac{\xi}{2NB} +\frac{r}{4N}\right)}=\delta_{k,\ell}\,I_M,\quad 0 \le k,\ell\le 2N-1.\eqno(4.16)$$

\parindent=0mm \vspace{.1in}
{\bf{Proof.}} Firstly, we will prove the necessary condition. By the orthonormality of the system $\{ {\bf \Psi}_k(t-\lambda)  e^{-\frac{- \iota \pi A}{B}(t^2-\lambda^2)}\}_{\lambda \in \Lambda,~k=0, 1,\dots,2N-1}$, we have

$$\Big\langle {\bf \Psi}_k(t-\lambda), {\bf \Psi}_{\ell}(t-\sigma)\Big\rangle= \int_{\mathbb R}{\bf \Psi}_k(t-\lambda)e^{-\frac{- \iota \pi A}{B}(t^2-\lambda^2)}\, \overline{{\bf \Psi}_{\ell}(t- \sigma) e^{-\frac{- \iota \pi A}{B}(t^2-\sigma^2)}} \,dt=e^{\iota \pi\frac{ A}{B}(\lambda^2-\sigma^2)}\delta_{k,\ell}\,\delta_{\lambda,\sigma}\,I_M,~  $$

\parindent=0mm \vspace{.1in}
where $\lambda,\sigma \in \Lambda$ and $k,\ell \in \left\{0, 1, 2,\dots, 2N-1\right\}$. Above relation can be recast in LCT domain as

$$\begin{array}{rcl}
\delta_{k,\ell}\, \delta_{\lambda, \sigma}I_{M} &=&\displaystyle \int_{\mathbb R} \hat{{\bf \Psi}}_k\left(\frac{\xi}{B}\right)e^{\frac{- 2 \pi \iota \xi \lambda}{B} }\,  \overline{\hat{{\bf \Psi}}_{\ell}\left(\frac{\xi}{B}\right)}e^{\frac{2 \pi \iota \xi \sigma}{B} } d\xi\\\\
&= &\displaystyle \int_{\mathbb R} \hat{{\bf \Psi}}_k\left(\frac{\xi}{B}\right)  \overline{\hat{{\bf \Psi}}_{\ell}\left(\frac{\xi}{B}\right)}e^{\frac{-2 \pi \iota \xi }{B}(\lambda-\sigma) } d\xi.
\end{array} $$

\parindent=0mm \vspace{.1in}
Taking $\lambda = 2m$ and $\sigma = 2n$, where $m, n \in \mathbb Z$, we have

$$\begin{array}{rcl}
\delta_{k,\ell}\, \delta_{m,n}I_M &=&\dfrac{1}{B}\displaystyle \int_{\mathbb R} e^{\frac{-4 \pi \iota \xi }{B}(m-n) }\hat{{\bf \Psi}}_k\left(\frac{\xi}{B}\right)  \overline{\hat{{\bf \Psi}}_{\ell}\left(\frac{\xi}{B}\right)}  d\xi \\\\
&=&\dfrac{1}{B}\displaystyle \int_{[0, BN]} e^{\frac{-4 \pi \iota \xi }{B}(m-n) }\sum_{j \in \mathbb Z}\hat{{\bf \Psi}}_k\left(\frac{\xi}{B}+Nj\right)  \overline{\hat{{\bf \Psi}}_{\ell}\left(\frac{\xi}{B}+Nj\right)}  d\xi
\end{array}$$

\parindent=0mm \vspace{.1in}
Define

$$F_{k,\ell}\left(\frac{\xi}{B}\right)=\sum_{j \in \mathbb Z}\hat{{\bf \Psi}}_k\left(\frac{\xi}{B}+Nj\right)  \overline{\hat{{\bf \Psi}}_{\ell}\left(\frac{\xi}{B}+Nj\right)} ,\quad 0\le k,\ell\le 2N-1. $$

\parindent=0mm \vspace{.1in}
Then, we have
$$\begin{array}{rcl}
\delta_{k,\ell}\,\delta_{m,n}I_M &=&\dfrac{1}{B}\displaystyle \int_{[0, BN]} e^{\frac{-4 \pi \iota \xi }{B}(m-n) } F_{k,\ell}\left(\frac{\xi}{B}\right)\, d\xi \\\\
&=&\dfrac{1}{B}\displaystyle \int_{[0, BN]} e^{\frac{-4 \pi \iota \xi }{B}(m-n) }\left\{ \sum_{s=0}^{2N-1}F_{k,\ell}\left(\frac{\xi}{B}+\frac{s}{2}\right)\right\}d\xi,
\end{array}$$

and
$$\sum_{s=0}^{2N-1}F_{k,\ell}\left(\frac{\xi}{B}+\frac{s}{2}\right)=2 \delta_{k,\ell}\,I_M. \qquad\eqno(4.17)$$

\parindent=0mm \vspace{.1in}
On taking $\lambda = \dfrac{r}{N}+2m$ and $\sigma = 2n$, where $m, n \in \mathbb Z$, we obtain

$$\begin{array}{rcl}
0&=&\displaystyle \int_{\mathbb R} \hat{{\bf \Psi}}_k\left(\frac{\xi}{B}\right)e^{\frac{- 2 \pi \iota \xi \lambda}{B} }\,  \overline{\hat{{\bf \Psi}}_{\ell}\left(\frac{\xi}{B}\right)}e^{\frac{2 \pi \iota \xi \sigma}{B} } d\xi \\\\
&=&\dfrac{1}{B}\displaystyle \int_{[0, BN]} e^{\frac{-2 \pi \iota \xi }{B}(\frac{r}{N}+2m+2n) }  \hat{{\bf \Psi}}_k\left(\frac{\xi}{B}\right)\overline{\hat{{\bf \Psi}}_{\ell}\left(\frac{\xi}{B}\right)} d\xi \\\\
&=&\dfrac{1}{B}\displaystyle \int_{[0, BN]} e^{-4 \pi \iota\frac{ \xi }{B}(m-n)} e^{-2 \pi \iota\frac{ \xi }{B}\frac{r}{N}} \sum_{j \in \mathbb Z}\hat{{\bf \Psi}}_k\left(\frac{\xi}{B}+Nj\right)  \overline{\hat{{\bf \Psi}}_{\ell}\left(\frac{\xi}{B}+Nj\right)}d\xi \\\\
&=&\dfrac{1}{B}\displaystyle \int_{[0, BN]} e^{-4 \pi \iota\frac{ \xi }{B}(m-n)} e^{-2 \pi \iota\frac{ \xi }{B}\frac{r}{N}} F_{k,\ell}\left(\frac{\xi}{B}\right)d\xi \\\\
&=&\dfrac{1}{B}\displaystyle \int_{[0, B/2)} e^{-4 \pi \iota\frac{ \xi }{B}(m-n)} e^{-2 \pi \iota\frac{ \xi }{B}\frac{r}{N}} \left\{\sum_{s=0}^{2N-1}e^{-2 \pi \iota\frac{r}{N}s} F_{k,\ell}\left(\frac{\xi}{B}+\frac{s}{2}\right)\right\}d\xi.
\end{array}$$

\parindent=0mm \vspace{.1in}
We conclude that

$$\sum_{s=0}^{2N-1}e^{-2 \pi \iota\frac{r}{N}s} F_{k,\ell}\left(\frac{\xi}{B}+\frac{s}{2}\right)={\bf 0}. \eqno(4.18)$$

\parindent=0mm \vspace{.1in}
Also we have

$$\sum_{j=0}^{2N-1}F_{k,\ell}\left(\frac{\xi}{B}+\frac{j}{2}\right)=\sum_{j \in \mathbb Z}\hat{{\bf \Psi}}_k\left(\frac{\xi}{B}+\frac{j}{2}\right)\overline{\hat{{\bf \Psi}}_{\ell}\left(\frac{\xi}{B}+\frac{j}{2}\right)}. $$

\parindent=0mm \vspace{.1in}
Therefore, equations (4.17) reduces to

$$\sum_{s\in \mathbb Z}\hat{{\bf \Psi}}_k\left(\frac{\xi}{B}+\frac{j}{2}\right)\overline{\hat{{\bf \Psi}}_{\ell}\left(\frac{\xi}{B}+\frac{j}{2}\right)}=2\delta_{k,\ell}\,I_M.\eqno(4.19)$$

\parindent=0mm \vspace{.1in}
Moreover, we have
$$\begin{array}{rcl}
F_{k,\ell}\left(\frac{2N\xi}{B}\right)&=&\displaystyle \sum_{j \in \mathbb Z} \hat{{\bf \Psi}}_k\left(2N\left(\frac{\xi}{B}+\frac{j}{2}\right) \right) \overline{\hat{{\bf \Psi}}_{\ell} \left(2N\left(\frac{\xi}{B}+\frac{j}{2}\right) \right)}\\\\
&=&\displaystyle\sum_{j \in \mathbb Z} H_k^\mu\left(\frac{\xi}{B}+\frac{j}{2}\right)\, \hat{{\bf \Phi}}\left(\frac{\xi}{B}+\frac{j}{2}\right)\overline{\hat{{\bf \Phi}}\left(\frac{\xi}{B}+\frac{j}{2}\right)}\,\overline{H^\mu_{\ell}\left(\frac{\xi}{B}+\frac{j}{2}\right)}\\\\
&=&\displaystyle\sum_{j =n.2N} H_k^\mu\left(\frac{\xi}{B}+nN\right)\, \hat{{\bf \Phi}}\left(\frac{\xi}{B}+nN\right)\overline{\hat{{\bf \Phi}}\left(\frac{\xi}{B}+nN\right)}\,\overline{H^\mu_{\ell}\left(\frac{\xi}{B}+nN\right)}\\\\
&&+\displaystyle\sum_{j =n.2N+1} H_k^\mu\left(\frac{\xi}{B}+nN+\frac{1}{2}\right)\, \hat{{\bf \Phi}}\left(\frac{\xi}{B}+nN+\frac{1}{2}\right)\\\
&&\qquad \qquad \qquad\overline{\hat{{\bf \Phi}}\left(\frac{\xi}{B}+nN+\frac{1}{2}\right)}\,\overline{H^\mu_{\ell}\left(\frac{\xi}{B}+nN+\frac{1}{2}\right)}\\\\
&&\quad +\cdots+\\\\
&&+\displaystyle\sum_{j =n.2N+(2N-1)} H_k^\mu\left(\frac{\xi}{B}+nN+\frac{2N-1}{2}\right)\, \hat{{\bf \Phi}}\left(\frac{\xi}{B}+nN+\frac{2N-1}{2}\right)\\\
&&\qquad \qquad \qquad\overline{\hat{{\bf \Phi}}\left(\frac{\xi}{B}+nN+\frac{2N-1}{2}\right)}\,\overline{H^\mu_{\ell}\left(\frac{\xi}{B}+nN+\frac{2N-1}{2}\right)}\\\\
&=&H_k^\mu\left(\frac{\xi}{B}\right)\left\{\displaystyle\sum_{j =n.2N}\hat{{\bf \Phi}}\left(\frac{\xi}{B}+nN\right)\overline{\hat{{\bf \Phi}}\left(\frac{\xi}{B}+nN\right)}\right\}\overline{H_{\ell}^\mu\left(\frac{\xi}{B}\right)}\\\
&&+H_k^\mu\left(\frac{\xi}{B}+\frac{1}{2}\right)\left\{\displaystyle\sum_{j =n.2N+1} \hat{{\bf \Phi}}\left(\frac{\xi}{B}+nN+\frac{1}{2}\right)\overline{\hat{{\bf \Phi}}\left(\frac{\xi}{B}+nN+\frac{1}{2}\right)}\right\}\overline{H_{\ell}^\mu\left(\frac{\xi}{B}+\frac{1}{2}\right)}\\\
&&+\cdots\\\
&&+H_k^\mu\left(\frac{\xi}{B}+\frac{2N-1}{2}\right)\left\{\displaystyle\sum_{j =n.2N+(2N-1)} \hat{{\bf \Phi}}\left(\frac{\xi}{B}+nN+\frac{2N-1}{2}\right)\overline{\hat{{\bf \Phi}}\left(\frac{\xi}{B}+nN+\frac{2N-1}{2}\right)}\right\}\\\
&&\qquad \qquad \qquad\overline{H_{\ell}^\mu\left(\frac{\xi}{B}+\frac{2N-1}{2}\right)}\\\\
&=&2\Big\{H_k^\mu\left(\frac{\xi}{B}\right)\overline{H_{\ell}^\mu\left(\frac{\xi}{B}\right)}+H_k^\mu\left(\frac{\xi}{B}+\frac{1}{2}\right)\overline{H_{\ell}^\mu\left(\frac{\xi}{B}+\frac{1}{2}\right)}+\cdots+\\\\
&&\qquad \qquad\qquad\qquad\qquad\qquad H_k^\mu\left(\frac{\xi}{B}+\frac{2N-1}{2}\right)\overline{H_{\ell}^\mu\left(\frac{\xi}{B}+\frac{2N-1}{2}\right)}\Big\}\\\
&=&2\displaystyle\sum_{j=0}^{2N-1}H_k^\mu\left(\frac{\xi}{B}+\frac{j}{2}\right)\overline{H_{\ell}^\mu\left(\frac{\xi}{B}+\frac{j}{2}\right)}.
\end{array}$$

\parindent=0mm \vspace{.1in}
Therefore, we have
$$\sum_{j \in \mathbb Z}\hat{{\bf \Psi}}_k\left(\frac{\xi}{B}+\frac{j}{2}\right)\overline{\hat{{\bf \Psi}}_{\ell}\left(\frac{\xi}{B}+\frac{j}{2}\right)}=2\displaystyle\sum_{j=0}^{2N-1}H_k^\mu\left(\frac{\xi}{2NB}+\frac{j}{4N}\right)\overline{H_{\ell}^\mu\left(\frac{\xi}{2NB}+\frac{j}{4N}\right)}.$$

\parindent=0mm \vspace{.1in}
By using (4.19), we conclude that

$$\sum_{j=0}^{2N-1}H_k^\mu\left(\frac{\xi}{2NB}+\frac{j}{4N}\right)\overline{H_{\ell}^\mu\left(\frac{\xi}{2NB}+\frac{j}{4N}\right)}=\delta_{k,\ell}\,I_M.$$

\parindent=8mm \vspace{.2in}
Now we will prove the sufficiency.

\parindent=0mm \vspace{.1in}
By equations (4.12), we have

$$\begin{array}{lcr}
\displaystyle\sum_{j \in \mathbb Z}\hat{{\bf \Psi}}_k\left(\frac{\xi}{B}+\frac{j}{2}\right)\overline{\hat{{\bf \Psi}}_{\ell}\left(\frac{\xi}{B}+\frac{j}{2}\right)}&&\\\\
\quad=\displaystyle\sum_{j\in\mathbb Z} H_k^\mu\left(\frac{\xi}{2NB}+\frac{j}{4N}\right)\hat{{\bf \Phi}}\left(\frac{\xi}{2NB}+\frac{j}{4N}\right)\overline{H_{\ell}^\mu\left(\frac{\xi}{2NB}+\frac{j}{4N}\right)} \, \overline{\hat{{\bf \Phi}}\left(\frac{\xi}{2NB}+\frac{j}{4N}\right)}&& \\\\
\quad =H_k^\mu\left(\frac{\xi}{2NB}+\frac{n}{2}\right)\left\{\displaystyle\sum_{j =n.2N} \hat{{\bf \Phi}}\left(\frac{\xi}{2NB}+\frac{n}{2}\right)\overline{\hat{{\bf \Phi}}\left(\frac{\xi}{2NB}+\frac{n}{2}\right)}\right\}\overline{H_{\ell}^\mu\left(\frac{\xi}{2NB}+\frac{n}{2}\right)} &&\\\\
\quad +H_k^\mu\left(\frac{\xi}{2NB}+\frac{n}{2}+\frac{1}{4N}\right)\left\{\displaystyle\sum_{j =n.2N} \hat{{\bf \Phi}}\left(\frac{\xi}{2NB}+\frac{n}{2}+\frac{1}{4N}\right)\overline{\hat{{\bf \Phi}}\left(\frac{\xi}{2NB}+\frac{n}{2}+\frac{1}{4N}\right)}\right\}\overline{H_{\ell}^\mu\left(\frac{\xi}{2NB}+\frac{n}{2}+\frac{1}{4N}\right)} &&\\\\
+\cdots&&\\\\
\quad +H_k^\mu\left(\frac{\xi}{2NB}+\frac{n}{2}+\frac{2N-1}{4N}\right)\left\{\displaystyle\sum_{j =n.2N} \hat{{\bf \Phi}}\left(\frac{\xi}{2NB}+\frac{n}{2}+\frac{2N-1}{4N}\right)\overline{\hat{{\bf \Phi}}\left(\frac{\xi}{2NB}+\frac{n}{2}+\frac{2N-1}{4N}\right)}\right\}&&\\\
\overline{H_{\ell}^\mu\left(\frac{\xi}{2NB}+\frac{n}{2}+\frac{2N-1}{4N}\right)} &&\\\\
\quad=2\,\left\{H_k^\mu\left(\frac{\xi}{2NB}\right)\overline{H_{\ell}^\mu\left(\frac{\xi}{2NB}\right)}+H_k^\mu\left(\frac{\xi}{2NB}+\frac{1}{4N}\right)\overline{H_{\ell}^\mu\left(\frac{\xi}{2NB}+\frac{1}{4N}\right)}+\cdots+\right.&&\\\\
\qquad\qquad\qquad \qquad\qquad\qquad \,\left. H_k^\mu\left(\frac{\xi}{2NB}+\frac{2N-1}{4N}\right)\overline{H_{\ell}^\mu\left(\frac{\xi}{2NB}+\frac{2N-1}{4N}\right)}\right\}&&\\\\
\quad=2\delta_{k,\ell}\,I_M.
\end{array}$$

\parindent=0mm \vspace{.1in}
It proves the orthonormality of the system $\big\{ {\bf \Psi}_k(x- \lambda)e^{-\frac{- \iota \pi A}{B}(t^2-\lambda^2)}:\lambda \in \Lambda, k=0, 1,\dots, 2N-1 \big\}.$
 \quad\fbox\\

\parindent=0mm \vspace{.2in}
{\bf Theorem 4.5.} {\it Suppose $\left\{ {\bf \Psi}_k(x- \lambda)e^{-\frac{- \iota \pi A}{B}(t^2-\lambda^2)} \right\}_{\lambda \in \Lambda,~k=0, 1,\dots, 2N-1}$ is the system as defined in Lemma 4.4 and orthonormal in $V_1$. Then this system is complete in $W_0^\mu\equiv V_1^\mu \ominus V_0^\mu$.}

\parindent=0mm \vspace{.2in}
{\bf Proof.} Since the system (4.15) is orthonormal in $V_1$. By Lemma 4.4 we have
$$\begin{array}{lcr},\left\{H_k^\mu\left(\frac{\xi}{2NB}\right)\overline{H_{\ell}^\mu\left(\frac{\xi}{2NB}\right)}+H_k^\mu\left(\frac{\xi}{2NB}+\frac{1}{4N}\right)\overline{H_{\ell}^\mu\left(\frac{\xi}{2NB}+\frac{1}{4N}\right)}+\cdots+\right.&&\\\\
\qquad\qquad\qquad \qquad\qquad\qquad \,\left. H_k^\mu\left(\frac{\xi}{2NB}+\frac{2N-1}{4N}\right)\overline{H_{\ell}^\mu\left(\frac{\xi}{2NB}+\frac{2N-1}{4N}\right)}\right\}&&\\\\
\quad=\delta_{k,\ell}\,I_M.
\end{array}$$

\parindent=0mm \vspace{.1in}
We will now prove its completeness.

\parindent=8mm \vspace{.2in}
For ${\bf f}_k \in W_0^\mu$, there exists constant matrices $\left\{P_{\lambda, k}^\mu\right\}$ such that

$${\bf f}_{k}(t)=\sqrt {2N}\sum_{\lambda \in \Lambda}P_{\lambda,k}^\mu\,{\bf \Phi}\big(2Nt-\lambda \big)e^{-\frac{- \iota \pi A}{B}(t^2-\lambda^2)},\quad 0\le k \le 2N-1.$$

\parindent=0mm \vspace{.1in}
Above relation can be written in the LCT domain as

$${\hat {\bf f}}_k\left(\frac{\xi}{B}\right)=P_k^\mu\left(\dfrac{\xi}{2NB}\right){\hat {\bf \Phi}}\left(\dfrac{\xi}{2NB}\right),\eqno(4.20)$$
where
$$P_k^\mu(\xi)=\dfrac{1}{\sqrt { qN}}\sum_{\lambda \in \Lambda}P_{\lambda,k}^\mu e^{-2 \pi \iota \lambda \xi /B}.$$

\parindent=0mm \vspace{.1in}
On the other hand, ${\bf f}_k\notin V_0^\mu$ and ${\bf f}_k\in W_0^\mu$ implies

$$\int_\mathbb R{\bf f}_k(t)\overline{{\bf \Phi}(t-\lambda)}e^{-\frac{- \iota \pi A}{B}(t^2-\lambda^2)}\,dt={\bf 0},\quad \lambda \in \Lambda.$$

\parindent=0mm \vspace{.1in}
This condition is equivalent to

$$\sum_{n \in \mathbb Z}{\hat {\bf f}}_k\left(\frac{\xi}{B}+\frac{n}{2}\right)\overline{{\hat{\bf \Phi}}\left(\frac{\xi}{B}+\frac{n}{2}\right)}={\bf 0},\quad \xi \in \mathbb R.$$

\parindent=0mm \vspace{.1in}
Therefore, the identities (4.7) and (4.20) give for all $\xi \in \mathbb R$,

$$\sum_{n \in \mathbb Z}P_k^\mu\left(\frac{\xi}{2NB}+\frac{j}{4N}\right)\hat{{\bf \Phi}}\left(\frac{\xi}{2NB}+\frac{j}{4N}\right)\overline{G^\mu\left(\frac{\xi}{2NB}+\frac{j}{4N}\right)} \, \overline{\hat{{\bf \Phi}}\left(\frac{\xi}{2NB}+\frac{j}{4N}\right)}={\bf 0}.$$

\parindent=0mm \vspace{.1in}
As similar to the identity (4.16) in Lemma 4.4, we have

$$P_k^\mu\left(\frac{\xi}{2NB}\right)\overline{G^\mu\left(\frac{\xi}{2NB}\right)}+P_k^\mu\left(\frac{\xi}{2NB}+\frac{1}{4N}\right)\overline{G^\mu\left(\frac{\xi}{2NB}+\frac{1}{4N}\right)}+\cdots+\qquad\qquad\qquad$$
$$P_k^\mu\left(\frac{\xi}{2NB}+\frac{2N-1}{4N}\right)\overline{G^\mu\left(\frac{\xi}{2NB}+\frac{2N-1}{4N}\right)}={\bf 0},\quad 0 \le k \le 2N-1.\eqno(4.21)$$

\parindent=0mm \vspace{.1in}
Let
$$\begin{array}{rcl}
\displaystyle P_{k^\prime}^\mu\left(\frac{\xi}{2NB}\right)&=&\displaystyle\overline{\left(P_k^\mu\left(\frac{\xi}{2NB}\right),P_k^\mu\left(\frac{\xi}{2NB}+\frac{1}{4N}\right),\dots,P_k^\mu\left(\frac{\xi}{2NB}+\frac{2N-1}{4N}\right)\right)},\\\\
\displaystyle \tilde G^\mu\left(\frac{\xi}{2NB}\right)&=&\displaystyle\overline{\left(G^\mu\left(\frac{\xi}{2NB}\right),G^\mu\left(\frac{\xi}{2NB}+\frac{1}{4N}\right),\dots,G^\mu\left(\frac{\xi}{2NB}+\frac{2N-1}{4N}\right)\right)},\\\\
\displaystyle H_{k^\prime}^\mu\left(\frac{\xi}{2NB}\right)&=&\displaystyle\overline{\left(H_k^\mu\left(\frac{\xi}{2NB}\right),H_k^\mu\left(\frac{\xi}{2NB}+\frac{1}{4N}\right),\dots,H_k^\mu\left(\frac{\xi}{2NB}+\frac{2N-1}{4N}\right)\right)}.
\end{array}$$

\parindent=0mm \vspace{.1in}
Then, equation (4.16) implies that for any $\xi \in \mathbb R$, the column vectors in $2NM\times M$ matrix $\tilde G^\mu$ and the column vectors in $2NM\times M$ matrix $H_{k^\prime}^\mu$ are orthogonal for $k=0,1,\dots,2N-1$ and these vectors form an orthogonal basis of $2NM$ dimensional complex Euclidean space $\mathbb C^{2NM}$.

\parindent=8mm \vspace{.1in}
Equation (4.21) implies that the column vectors in $2NM\times M$ matrix $P_{k^\prime}^\mu$ and the column vectors of $2NM\times M$ matrix $\tilde G^\mu$ are orthogonal. Therefore, there exists an $M \times M$ matrix $Q_k(\xi)$ such that

$$P_k^\mu\left(\frac{\xi}{B}\right)=Q_k^\mu\left(\frac{\xi}{B}\right)H_k^\mu\left(\frac{\xi}{B}\right),\quad \xi \in \mathbb R,~ 0\le k \le 2N-1.$$

\parindent=0mm \vspace{.2in}
Therefore, from equations (4.12) and (4.20), we have
$$\begin{array}{rcl}
{\hat {\bf f}}_k\left(\frac{\xi}{B}\right)&=&P_k^\mu\left(\frac{\xi}{2NB}\right){\hat {\bf \Phi}}\left(\frac{\xi}{2NB}\right)\\\\
&=&Q_k^\mu\left(\frac{\xi}{2NB}\right)H_k^\mu\left(\frac{\xi}{2NB}\right){\hat {\bf \Phi}}\left(\frac{\xi}{2NB}\right)\\\\
&=&Q_k^\mu\left(\frac{\xi}{2NB}\right){\hat {\bf \Psi}}_k\left(\frac{\xi}{B}\right).
\end{array}$$

\parindent=0mm \vspace{.1in}
By using the orthonormality of the system (4.15), we have

$$\int_\mathbb R{\hat {\bf f}}_k\left(\frac{2N\xi}{B}\right)\overline{{\hat {\bf f}}_k\left(\frac{2N\xi}{B}\right)}\,d\xi=\int_{\mathbb R} Q_k^\mu\left(\frac{\xi}{B}\right){\hat {\bf \Psi}}_k\left(\frac{2N\xi}{B}\right)\overline{{\hat {\bf \Psi}}_k\left(\frac{2N\xi}{B}\right)}\,\overline{Q_k^\mu\left(\frac{\xi}{B}\right)}\,d\xi.$$

\parindent=0mm \vspace{.1in}
Therefore, we have
$$\int_\mathbb R{\hat {\bf f}}_k\left(\frac{2N\xi}{B}\right)\overline{{\hat {\bf f}}_k\left(\frac{2N\xi}{B}\right)}\,d\xi=2\int_0^{1/2}Q_k^\mu\left(\frac{\xi}{B}\right)\overline{Q_k^\mu\left(\frac{\xi}{B}\right)}\,d\xi.$$

\parindent=0mm \vspace{.1in}
This shows that $P_k(^\mu\xi)$ has the series expansion and let the constant $M \times M$ matrices $\left\{R_{\lambda,k}^\mu\right\}_{\lambda \in \Lambda,k=0,1,\dots,2N-1}$ be its coefficients. Therefore, we have

$${\bf f}_k(t)=\sum_{\lambda \in \Lambda}R_{\lambda,k}^\mu\,{\bf \Psi}_k(t-\lambda)e^{-\frac{- \iota \pi A}{B}(t^2-\lambda^2)}.$$

\parindent=0mm \vspace{.1in}
It proves the completeness of the system $\left\{ {\bf \Psi}_k(x- \lambda) e^{-\frac{- \iota \pi A}{B}(t^2-\lambda^2)}\right\}_{\lambda \in \Lambda,~k=0, 1,\dots, 2N-1}$ in $W_0$. \qquad\fbox

\parindent=8mm\vspace{.2in}
If ${\bf \Psi}_0^\mu,{\bf \Psi}_1^\mu,\dots,{\bf \Psi}_{2 N-1}^\mu \in V_1^\mu  $ are as in Lemma 4.4, one can obtain from them as orthonormal basis for $L^2\big(\mathbb R, \mathbb C^{M}\big)$ by following the standard procedure for construction of wavelet from a given MRA. It can be easily checked that for every $j\in \mathbb{Z}$, the collection $\left\{ \sqrt {2N} {\bf \Psi}_k \left(\left(2N \right)^j t - \lambda  \right)e^{-\frac{- \iota \pi A}{B}(t^2-\lambda^2)}:\lambda \in \Lambda, k=0, 1,\dots, 2 N-1  \right\}$
is a complete orthogonal system for $V_{j+1}^\mu$. Therefore, it follows immediately from (4.5) that the collection

$$\left\{ \sqrt {2N} {\bf \Psi}_k \left(\left(2N \right)^j t - \lambda  \right)e^{-\frac{- \iota \pi A}{B}(t^2-\lambda^2)}:\lambda \in \Lambda, k=0, 1,\dots, 2 N-1  \right\}$$

\parindent=0mm \vspace{.1in}
forms a complete orthonormal system for $L^2\big(\mathbb R, \mathbb C^{M}\big)$.

\parindent=0mm \vspace{.2in}
{\bf{5. Construction of LCT-VNUMRA }}

\parindent=0mm \vspace{.2in}
The main goal of this section is to construct a LCT-VNUMRA starting from a vector-valued refinement mask $G(\xi)$ of the form

$$G^\mu\left(\frac{\xi}{B}\right)=G_{\lambda, 1}^\mu\left(\frac{\xi}{B}\right)+e^{-2 \pi \iota \frac{r}{N}\frac{\xi}{B}}G_{\lambda,2}^\mu\left(\frac{\xi}{B}\right), \eqno(5.1)$$

\parindent=0mm \vspace{.1in}
where $N >1$ is an integer and $r$ is an odd integer with $1 \le r \le 2N -1$ such that $r$ and $N$ are relatively prime and $G_{\lambda, 1}^\mu\left(\frac{\xi}{B}\right)$
 and $G_{\lambda, 2}^\mu\left(\frac{\xi}{B}\right)$ are $M\times M$ constant symmetric matrix sequences. In other words, we establish conditions under which the solutions of scaling equation (4.6) generate a LCT-VNUMRA in $L^2(\mathbb R)$ or equivalently, we find a sufficient for the orthonormality of the system $\big\{{\bf \Phi}(t-\lambda)e^{-\frac{- \iota \pi A}{B}(t^2-\lambda^2)}:\lambda \in \Lambda\big\} $, where $\Lambda=\left\{0, r/N\right\}+2 \mathbb Z$. Therefore, the scaling vector $\Phi$ associated with given LCT-VNUMRA should satisfy the scaling identity

$$\hat{\bf \Phi}\left(\frac{2N\xi}{B}\right)=G^\mu\left(\frac{\xi}{B}\right)\hat{\bf \Phi}\left(\frac{\xi}{B}\right). \eqno(5.2)$$

\parindent=0mm \vspace{.1in}
We further assume that:

$$\sum_{s=0}^{2N-1}G^\mu\left(\frac{\xi}{2NB}+\frac{s}{4N}\right)\overline{G^\mu\left(\frac{\xi}{2NB}+\frac{s}{4N}\right)}=I_M. \eqno(5.3)$$

\parindent=0mm \vspace{.2in}
{\bf Theorem 5.1.} {\it Let $G^\mu\left(\frac{\xi}{B}\right)$ be the vector-valued refinement mask associated with the  vector-valued scaling function ${\bf \Phi}$ of LCT-VNUMRA and satisfies the condition (5.3) together with $G^\mu(0)=I_M$ and\, $G^\mu\left(\frac{\xi}{B}\right)=\overline{G^\mu\left(\frac{\xi}{B}\right)}, \forall~\xi\in \mathbb R$. Then, a sufficient condition  for the collection $\left\{{\bf \Phi}(x-\lambda)e^{-\frac{- \iota \pi A}{B}(t^2-\lambda^2)}\right\}_{: \lambda\in \Lambda}$ to be orthonormal in $L^2\big(\mathbb R, \mathbb C^{M}\big)$ is the existence of a constant $C>0$ and of a compact set $E\subset \mathbb R$ that contains the neighbourhood of the origin such that

$$\left|G^\mu\left(\frac{\xi}{(2N)^kB}\right)\right| \ge C,\quad \forall~\xi\in \mathbb R, k\in\mathbb Z.\eqno(5.4)$$ }

\parindent=0mm \vspace{.2in}
{\bf Proof.} Let us assume the existence of a constant $C$ and of the compact set $E\subset K$ with properties satisfied above. For any $k\in \mathbb N$, we define

$$g_k\left(\frac{\xi}{B}\right)=\left\{\prod_{j=1}^{k} G^\mu\left(\frac{\xi}{(2N)^kB}\right)\right\}\chi_{E}\left(\frac{\xi}{(2N)^kB}\right).$$

\parindent=0mm \vspace{.1in}
As the interior of the compact set $E$ contains ${\bf 0}, g_k \to \hat{{\bf \Phi}}$ pointwise as $k \to \infty $. Therefore, there exists a constant $W>0$ such that $\big|G^\mu\left(\frac{\xi}{B}\right)-G^\mu(0)\big|\le W|\xi|,\,\text{for all}\,\xi \in \mathbb R,$ and thus $\big|G^\mu\left(\frac{\xi}{B}\right)\big|\ge 1-W|\xi|.$ Since $E$ is bounded, we can find an integer $k_{0}\in \mathbb Z$  such that $W|\xi|< (2N)^{k}, \text{for}\,k>k_{0}, \xi\in E$ and hence, there exists a constant $C_{1}>0$ such that

$$\chi_{E}\left(\frac{\xi}{B}\right)\le C_{1}\big|\hat{{\bf \Phi}}\left(\frac{\xi}{B}\right)\big|, \quad \text{for all}~\xi \in \mathbb R.$$

Thus, we have
$$\left|g_{k}\left(\frac{\xi}{B}\right)\right|\le C_{1}\left\{\prod_{j=1}^{k}\left| G^\mu\left(\frac{\xi}{(2N)^kB}\right)\right|\right\}\left|\hat{{\bf \Phi}}\left(\frac{\xi}{(2N)^kB}\right)\right|=C_{1}\left|\hat{{\bf \Phi}}\left(\frac{\xi}{B}\right)\right|.$$

\parindent=0mm \vspace{.1in}
Therefore, by Lebesgue dominated convergence theorem the sequence $\left\{g_{k}\right\}$ converges to $\hat{\bf \Phi}$ in $L^2$-norm. We will now compute by induction the integral

$$\int_\mathbb Rg_{k}(\xi)\,\overline{g_{k}(\xi)}\,\overline{\chi_{(\lambda-\sigma)}(\xi)}\,d\xi,\quad \text{where}~\lambda,\sigma \in \Lambda.$$

\parindent=0mm \vspace{.1in}
For $k=1$, we have
$$\begin{array}{lcr}
\displaystyle\int_\mathbb R g_{1}\left(\frac{\xi}{B}\right)\,\overline{g_{1}\left(\frac{\xi}{B}\right)}e^{-2 \pi \iota\frac{ \xi}{B}(\lambda -\sigma)}\,d\xi\\\\
\qquad=\displaystyle\int_\mathbb R G^\mu\left(\frac{\xi}{2NB}\right)\overline{G^\mu\left(\frac{\xi}{2NB}\right)}\,{\chi_E\left(\frac{\xi}{2NB}\right)}e^{-2 \pi \iota\frac{ \xi}{B}(\lambda -\sigma)}\,d\xi\\\\
\qquad=\displaystyle(2N)\displaystyle\int_E G^\mu\left(\frac{\xi}{B}\right)\, \overline{G^\mu\left(\frac{\xi}{B}\right)}e^{-2 \pi \iota\frac{(2N) \xi}{B}(\lambda -\sigma)}\,d\xi\\\\
\qquad=4N\displaystyle\int_0^{B/2}\left\{\sum_{s=0}^{2N-1} G^\mu\left(\frac{\xi}{B}+\frac{1}{4N}\right)\overline{G^\mu\left(\frac{\xi}{B}+\frac{1}{4N}\right)}\,e^{-\pi \iota\frac{\xi}{B}(\lambda -\sigma)s}\right\}\\\
\qquad\qquad\qquad\qquad\qquad\qquad\qquad\qquad\qquad\qquad\qquad\qquad\qquad\times~e^{-2 \pi \iota\frac{(2N) \xi}{B}(\lambda -\sigma)}\,d\xi.
\end{array}$$

\parindent=0mm \vspace{.1in}
If $\lambda-\sigma \in2 \mathbb Z$, then the expression in the brackets in the above integral is equal to $I_{M}$ by (5.3) and thus

$$\begin{array}{rcl}
\displaystyle\int_\mathbb R g_{1}\left(\frac{\xi}{B}\right)\,\overline{g_{1}\left(\frac{\xi}{B}\right)}e^{-2 \pi \iota\frac{ \xi}{B}(\lambda -\sigma)}\,d\xi&=&\displaystyle 4N\displaystyle\int_0^{B/4N}I_M e^{-2 \pi \iota\frac{(2N) \xi}{B}(\lambda -\sigma)} d\xi \\\\
&=&\displaystyle 2\int_{[0, B/2)}I_M e^{-2 \pi \iota\frac{ \xi}{B}(\lambda -\sigma)} d\xi\\\\
&=&\displaystyle\delta_{\lambda,\sigma}\,I_M.
\end{array}$$

\parindent=0mm \vspace{.1in}
On the other hand, if $\lambda=2m,\, \sigma=2n+r/N,$ where $m,n \in \mathbb Z$, then the same expression will vanish and the integral becomes

$$\int_\mathbb R g_{1}\left(\frac{\xi}{B}\right)\,\overline{g_{1}\left(\frac{\xi}{B}\right)}e^{-2 \pi \iota\frac{ \xi}{B}(\lambda -\sigma)}\,d\xi={\bf 0}.$$

\parindent=0mm \vspace{.1in}
When $k\ge 2$, we have

$$\begin{array}{lrc}
\displaystyle\int_\mathbb R g_{k}\left(\frac{\xi}{B}\right)\,\overline{g_{k}\left(\frac{\xi}{B}\right)}e^{-2 \pi \iota\frac{ \xi}{B}(\lambda -\sigma)}\,d\xi\\\\
\quad=\displaystyle\int_\mathbb R G^\mu\left(\frac{\xi}{(2N)^1B}\right)G^\mu\left(\frac{\xi}{(2N)^2B}\right)\dots G^\mu\left(\frac{\xi}{(2N)^kB}\right)\overline{G^\mu\left(\frac{\xi}{(2N)^kB}\right)}\,\overline{G^\mu\left(\frac{\xi}{(2N)^{k-1}B}\right)}\\\\
\quad\qquad\qquad\qquad\qquad\qquad\qquad\qquad\qquad\qquad\dots~\overline{G^\mu\left(\frac{\xi}{(2N)^1B}\right)}e^{-2 \pi \iota\frac{ \xi}{B}(\lambda -\sigma)}d\xi \\\\
\quad=(2N)^{k}\displaystyle\int_E G^\mu\left(\frac{(2N)^{k-1}\xi}{B}\right)\,G^\mu\left(\frac{(2N)^{k-2}\xi}{B}\right)\dots G^\mu\left(\frac{\xi}{B}\right)\,\overline{G^\mu\left(\frac{\xi}{B}\right)}\,\overline{G^\mu\left(\frac{(2N)\xi}{B}\right)}\\\
\quad\qquad\qquad\qquad\qquad\qquad\qquad\qquad\qquad\qquad\qquad\dots~\overline{G^\mu\left(\frac{(2N)^{k-1}\xi}{B}\right)}e^{-2 \pi \iota\frac{ (2N)^k\xi}{B}(\lambda -\sigma)}d\xi\\\\
\quad =(2N)^k\displaystyle\int_E \left\{\displaystyle\prod_{\ell=0}^{k-1}G^\mu\left(\frac{(2N)^\ell\xi}{B}\right)\right\}\overline{\left\{\displaystyle\prod_{\ell=0}^{k-1}G^\mu\left(\frac{(2N)^\ell\xi}{B}\right)\right\}}e^{-2 \pi \iota\frac{ (2N)^k\xi}{B}(\lambda -\sigma)}d\xi\\\\
\quad =(2N)^k\displaystyle\int_E \left\{\displaystyle\prod_{\ell=1}^{k-1}G^\mu\left(\frac{(2N)^\ell\xi}{B}\right)\right\}G^\mu\left(\frac{\xi}{B}\right)\,\overline{G^\mu\left(\frac{\xi}{B}\right)}\overline{\left\{\displaystyle\prod_{\ell=1}^{k-1}G^\mu\left(\frac{(2N)^\ell\xi}{B}\right)\right\}}e^{-2 \pi \iota\frac{ (2N)^k\xi}{B}(\lambda -\sigma)}d\xi\\\\
\quad =2(2N)^k\displaystyle\int_0^{B/2} \left\{\displaystyle\prod_{\ell=1}^{k-1}G^\mu\left(\frac{(2N)^\ell\xi}{B}\right)\right\}G^\mu\left(\frac{\xi}{B}\right)\,\overline{G^\mu\left(\frac{\xi}{B}\right)}\overline{\left\{\displaystyle\prod_{\ell=1}^{k-1}G^\mu\left(\frac{(2N)^\ell\xi}{B}\right)\right\}}e^{-2 \pi \iota\frac{ (2N)^k\xi}{B}(\lambda -\sigma)}d\xi\\\\
\quad =2(2N)^k\displaystyle\int_0^{B/4N} \left\{\displaystyle\prod_{\ell=2}^{k-1}G^\mu\left(\frac{(2N)^\ell\xi}{B}\right)\right\}\left[\sum_{s=0}^{2N-1}G^\mu\left(\frac{2N\xi}{B}+\frac{s}{2}\right)\,G^\mu\left(\frac{\xi}{B}+\frac{s}{2}\right)\right.\\\\
\qquad\qquad\qquad\left.\overline{G^\mu\left(\frac{2N\xi}{B}+\frac{s}{2}\right)\;G^\mu\left(\frac{\xi}{B}+\frac{s}{2}\right)}\right]\overline{\left\{\displaystyle\prod_{\ell=1}^{k-1}G^\mu\left(\frac{(2N)^\ell\xi}{B}\right)\right\}}e^{-2 \pi \iota\frac{ (2N)^k\xi}{B}(\lambda -\sigma)}d\xi\\\\
\quad =2(2N)^k\displaystyle\int_0^{B/4N} \left\{\displaystyle\prod_{\ell=2}^{k-1}G^\mu\left(\frac{(2N)^\ell\xi}{B}\right)\right\} S\left(\frac{\xi}{B}\right) \overline{\left\{\displaystyle\prod_{\ell=2}^{k-1}G^\mu\left(\frac{(2N)^\ell\xi}{B}\right)\right\}}e^{-2 \pi \iota\frac{ (2N)^k\xi}{B}(\lambda -\sigma)}d\xi
\end{array}$$

\parindent=0mm \vspace{.1in}
where
$$S\left(\frac{\xi}{B}\right)=\left\{\sum_{s=0}^{2N-1}G^\mu\left(\frac{2N\xi}{B}+\frac{s}{2}\right)\,G^\mu\left(\frac{\xi}{B}+\frac{s}{2}\right)\overline{G^\mu\left(\frac{2N\xi}{B}+\frac{s}{2}\right)\;G^\mu\left(\frac{\xi}{B}+\frac{s}{2}\right)}\right\}.$$

\parindent=0mm \vspace{.1in}
Since the refinement mask $G^\mu\left(\frac{\xi}{B}\right)$ can be expressed as (5.1), therefore, the above relation becomes

$$\begin{array}{rcl}
S\left(\frac{\xi}{B}\right)&=&G_{\lambda,1}^\mu\left(\frac{\xi}{B}\right)\,\overline{G_{\lambda,1}^\mu\left(\frac{\xi}{B}\right)}\,G_{\lambda,2}^\mu\left(\frac{\xi}{B}\right)\,\overline{G_{\lambda,2}^\mu\left(\frac{\xi}{B}\right)}\\\
&=&G^\mu\left(\frac{2N\xi}{B}\right)\overline{G^\mu\left(\frac{2N\xi}{B}\right)}.
\end{array}$$

\parindent=0mm \vspace{.1in}
Thus, we have
$$\begin{array}{lcr}
\displaystyle\int_\mathbb R g_{k}\left(\frac{\xi}{B}\right)\,\overline{g_{k}\left(\frac{\xi}{B}\right)}e^{-2 \pi \iota\frac{ \xi}{B}(\lambda -\sigma)}\,d\xi\\\\
\qquad =(2N)^k\displaystyle\int_{E/2N} \left\{\displaystyle\prod_{\ell=1}^{k-1}G^\mu\left(\frac{(2N)^\ell\xi}{B}\right)\right\} \overline{\left\{\displaystyle\prod_{\ell=1}^{k-1}G^\mu\left(\frac{(2N)^\ell\xi}{B}\right)\right\}}e^{-2 \pi \iota\frac{ (2N)^k\xi}{B}(\lambda -\sigma)}d\xi\\\
\qquad =(2N)^{k-1}\displaystyle\int_{E} \left\{\displaystyle\prod_{\ell=0}^{k-2}G^\mu\left(\frac{(2N)^\ell\xi}{B}\right)\right\} \overline{\left\{\displaystyle\prod_{\ell=0}^{k-2}G^\mu\left(\frac{(2N)^\ell\xi}{B}\right)\right\}}e^{-2 \pi \iota\frac{ (2N)^k\xi}{B}(\lambda -\sigma)}d\xi\\\
\qquad =\displaystyle\int_\mathbb R g_{k-1}\left(\frac{\xi}{B}\right)\,\overline{g_{k-1}\left(\frac{\xi}{B}\right)}e^{-2 \pi \iota\frac{ \xi}{B}(\lambda -\sigma)}\,d\xi
\end{array}$$

\parindent=0mm \vspace{.1in}
Therefore for any $k\in\mathbb Z$, we have

$$\displaystyle\int_\mathbb R g_{k}\left(\frac{\xi}{B}\right)\,\overline{g_{k}\left(\frac{\xi}{B}\right)}e^{-2 \pi \iota\frac{ \xi}{B}(\lambda -\sigma)}\,d\xi=\delta_{\lambda, \sigma},\quad \lambda, \sigma \in \Lambda.$$

\parindent=0mm \vspace{.1in}
Passing to the limit as $k\to \infty$ and using Plancherel's formula, we obtain

$$\int_\mathbb R{\bf \Phi}(x-\lambda)e^{\frac{-2 \iota \pi A}{B}(t^2-\lambda^2)}\overline{{\bf \Phi}(x-\sigma)e^{\frac{- 2\iota \pi A}{B}(t^2-\sigma^2)}}\,dx=\int_\mathbb R{\hat {\bf \Phi}}(\xi)\overline{{\hat {\bf \Phi}}(\xi)}e^{\frac{-2 \iota \pi A}{B}(\lambda^2-\sigma^2)}\,d\xi=\delta_{\lambda, \sigma},\quad \lambda, \sigma \in \Lambda$$

\parindent=0mm \vspace{.1in}
which proves the desired orthonormality. \qquad\qquad \fbox

\parindent=0mm \vspace{.2in}
{\bf {Acknowledgments }}

\parindent=0mm \vspace{.2in}
This work  is supported by the UGC-BSR Research Start Up Grant(No. F.30-498/2019(BSR)) provided by UGC, Govt. of India.

\parindent=0mm \vspace{.2in}
{\bf{References}}

\begin{enumerate}

{\small {
\item Abdullah, Vector-valued multiresolution analysis on local fields, {\it Analysis.} {\bf 34}(4) (2014)  415-428.

\item B. Behera and Q. Jahan, Multiresolution analysis on local fields and characterization of scaling functions, {\it Adv. Pure Appl. Math.} {\bf 3} (2012) 181-202.

\item   J. J. Benedetto and R. L. Benedetto, Wavelet theory for local fields and related groups, {\it J. Geom. Anal.} {\bf 14} (2004) 424-456.

\item  Q. Chen and Z. Chang, A study on compactly supported orthogonal vector-valued wavelets and wavelet packets, {\it Chaos, Solitons Fract.} {\bf 31} (2007) 1024-1034.

\item L. Debnath and F. A. Shah, {\it Wavelet Transforms and Their Applications}, Birkh\"{a}user, New York, 2015.

\item Yu. A. Farkov,  Orthogonal wavelets with compact support on locally compact Abelian groups,  {\it Izv. Math.} {\bf 69}(3) (2005) 623-650.

\item  J. P. Gabardo and M. Nashed, Nonuniform multiresolution analyses and spectral pairs, {\it J. Funct. Anal.} {\bf 158} (1998)  209-241.

\item  J. P. Gabardo and M. Nashed, An analogue of Cohen's condition for nonuniform multiresolution analyses, in:  {\it Wavelets, Multiwavelets and Their Applications}, A. Aldroubi, E. Lin (Eds.), Amer. Math. Soc., Providence, RI, (1198) 41-61.

\item  H. K. Jiang, D. F. Li and N. Jin, Multiresolution analysis on local fields, {\it J. Math. Anal. Appl.} {\bf 294} (2204) 523-532.

\item A. Yu. Khrennikov,V. M. Shelkovich and M. Skopina, $p$-Adic refinable functions and MRA-based wavelets, {\it J. Approx. Theory,}  {\bf 161} (2009)  226-238.

\item W. C. Lang, Orthogonal wavelets on the Cantor dyadic group, {\it SIAM J. Math. Anal.} {\bf 27} (1996) 305-312.

\item S. F. Lukomskii, Step refinable functions and orthogonal MRA on Vilenkin groups.  {\it J. Fourier Anal. Appl.} {\bf 20} (2014) 42-65.

\item S. G. Mallat, Multiresolution approximations and wavelet orthonormal bases of $L^2(\mathbb R)$,  {\it Trans. Amer. Math. Soc.} {\bf 315} (1989)  69-87.

\item Meenakshi, P. Manchanda and A. H. Siddiqi,  Wavelets associated with nonuniform multiresolution analysis on positive half-line, {\it Int. J. Wavelets Multiresolut. Inf. Process.} {\bf 10}(2) (2012) 1250018, 27pp.

\item G. \'{O}lafsson, Continuous action of Lie groups on $\mathbb R^n$ and frames, {\it Int. J. Wavelets Multiresolut. Inf. Process.} {\bf 3} (2005) 211-232.

\item D. Ramakrishnan and R. J. Valenza, {\it Fourier Analysis on Number Fields}, Graduate Texts in Mathematics 186 (Springer-Verlag, New York, 1999)

\item F. A. Shah, Construction of wavelet packets on $p$-adic field, {\it Int. J. Wavelets Multiresolut. Inf. Process.} {\bf 7}(5) (2009)  553-565.

\item F. A. Shah, Tight wavelet frames generated by the Walsh polynomials, {\it Int. J. Wavelets Multiresolut. Inf. Process.} {\bf 11}(6) (2013), pp. 15 pages.

\item F. A. Shah, Frame multiresolution analysis on local fields of positive characteristic, {\it J. Operators}. Article ID 216060, 8 pages (2015).

\item F. A. Shah and Waseem, Nonuniform Multiresolution Analysis Associated with Linear Canonical Transform ,  {\it preprint.}  (2020).

\item F. A. Shah and Abdullah, A characterization of tight wavelet frames on local fields of positive characteristic, {\it J. Contemp. Math. Anal.} {\bf 49} (2014) 251-259.

\item F. A. Shah and Abdullah, Nonuniform multiresolution analysis on local fields of positive characteristic,  {\it Comp. Anal. Opert. Theory.}  (2014)  DOI 10.1007/s11785-014-0412-0.

\item F. A. Shah and L. Debnath, Tight wavelet frames on local fields, {\it Analysis}, {\bf 33} (2013)  293-307.

\item  M. H. Taibleson, {\it Fourier Analysis on Local Fields} (Princeton University Press, Princeton, 1975).

\item X. G. Xia and B. W. Suter, Vector-valued wavelets and vector filter banks, {\it IEEE Trans. Signal Process.} {\bf 44}(3) (1996) 508-518.

\item P. Z. Xie and J. X. He, Multiresolution analysis and Haar wavelets on the product of Heisenberg group, {\it Int. J. Wavelets Multiresolut. Inf. Process.} {\bf 7}(2) (2009) 243-254.
}}

\end{enumerate}

\end {document}